\documentclass[12pt]{amsart}      
\usepackage{txfonts}
\usepackage{amssymb}
\usepackage{eucal}
\usepackage{graphicx}
\usepackage{amsmath}
\usepackage{amscd}
\usepackage[all]{xy}           
\usepackage[active]{srcltx} 

\usepackage{amsfonts,latexsym}
\usepackage{xspace}
\usepackage{epsfig}
\usepackage{float}
\usepackage{color}
\usepackage{fancybox}
\usepackage{colordvi}
\usepackage{multicol}
\usepackage{colordvi}
\usepackage[hypertex]{hyperref} 

\newtheorem{theorem}{Theorem}[section]
\newtheorem{definition}[theorem]{Definition}
\newtheorem{lemma}[theorem]{Lemma}

\newtheorem{prop-def}{Proposition-Definition}[section]
\newtheorem{coro-def}{Corollary-Definition}[section]

\newtheorem{example}{Example}[section]

\newtheorem{problem}[theorem]{Problem}

\newcommand{\nc}{\newcommand}
\nc{\tred}[1]{\textcolor{red}{#1}}
\nc{\tblue}[1]{\textcolor{blue}{#1}}
\nc{\tgreen}[1]{\textcolor{green}{#1}}
\nc{\tpurple}[1]{\textcolor{purple}{#1}}
\nc{\btred}[1]{\textcolor{red}{\bf #1}}
\nc{\btblue}[1]{\textcolor{blue}{\bf #1}}
\nc{\btgreen}[1]{\textcolor{green}{\bf #1}}
\nc{\btpurple}[1]{\textcolor{purple}{\bf #1}}

\renewcommand{\Bbb}{\mathbb}

\newcommand{\efootnote}[1]{}

\renewcommand{\textbf}[1]{}

\newcommand{\delete}[1]{}
\nc{\dfootnote}[1]{{}}          
\nc{\ffootnote}[1]{\dfootnote{#1}}
\nc{\mfootnote}[1]{\footnote{#1}} 
\nc{\ofootnote}[1]{\footnote{\tiny Older version: #1}}

\nc{\mlabel}[1]{\label{#1}}  
\nc{\mcite}[1]{\cite{#1}}  
\nc{\mref}[1]{\ref{#1}}  
\nc{\mbibitem}[1]{\bibitem{#1}} 

\delete{
\nc{\mlabel}[1]{\label{#1}  
{\hfill \hspace{1cm}{\bf{{\ }\hfill(#1)}}}}
\nc{\mcite}[1]{\cite{#1}{{\bf{{\ }(#1)}}}}  
\nc{\mref}[1]{\ref{#1}{{\bf{{\ }(#1)}}}}  
\nc{\mbibitem}[1]{\bibitem[\bf #1]{#1}} 
}

\nc{\mtail}{\leq_t}
\nc{\mhead}{\leq_h}
\nc{\rk}{\mathrm{rk}}
\nc{\mset}[1]{\tilde{#1}}
\nc{\pa}{\frakL}
\nc{\arr}{\rightarrow}
\nc{\lu}[1]{(#1)}
\nc{\mult}{\mrm{mult}}
\nc{\diff}{\mathrm{Der}}
\nc{\indiff}{\mathrm{InDer}}
\nc{\outdiff}{\mathrm{OutDer}}
\nc{\conmat}{connection matrix\xspace}
\nc{\bounmat}{boundary matrix\xspace}
\nc{\pcyc}{\mathfrak c}
\nc{\calpa}{\calp_A}
\nc{\calpal}{\Gamma_{AL}}
\nc{\calpc}{\calp_L}
\nc{\frakDa}{\frakD_1}
\nc{\frakDal}{\frakD_2}
\nc{\frakDc}{\frakD_L}
\nc{\frakDv}{\frakD_V}
\nc{\frakDp}{\frakD_F}
\nc{\frakBa}{\frakB_1}
\nc{\frakBal}{\frakB_2}
\nc{\frakBc}{\frakB_L}
\nc{\frakBv}{\frakB_V}

\nc{\bin}[2]{ (_{\stackrel{\scs{#1}}{\scs{#2}}})}  
\nc{\binc}[2]{ \left (\!\! \begin{array}{c} \scs{#1}\\
    \scs{#2} \end{array}\!\! \right )}  
\nc{\bincc}[2]{  \left ( {\scs{#1} \atop
    \vspace{-1cm}\scs{#2}} \right )}  
\nc{\bs}{\bar{S}}
\nc{\cosum}{\sqsubset}
\nc{\la}{\longrightarrow}
\nc{\rar}{\rightarrow}
\nc{\dar}{\downarrow}
\nc{\dprod}{**}
\nc{\dap}[1]{\downarrow \rlap{$\scriptstyle{#1}$}}
\nc{\md}{\mathrm{dth}}
\nc{\uap}[1]{\uparrow \rlap{$\scriptstyle{#1}$}}
\nc{\defeq}{\stackrel{\rm def}{=}}
\nc{\disp}[1]{\displaystyle{#1}}
\nc{\dotcup}{\ \displaystyle{\bigcup^\bullet}\ }
\nc{\gzeta}{\bar{\zeta}}
\nc{\hcm}{\ \hat{,}\ }
\nc{\hts}{\hat{\otimes}}
\nc{\barot}{{\otimes}}
\nc{\free}[1]{\bar{#1}}
\nc{\uni}[1]{\tilde{#1}}
\nc{\hcirc}{\hat{\circ}}
\nc{\lleft}{[}
\nc{\lright}{]}
\nc{\lc}{\lfloor}
\nc{\rc}{\rfloor}
\nc{\curlyl}{\left \{ \begin{array}{c} {} \\ {} \end{array}
    \right .  \!\!\!\!\!\!\!}
\nc{\curlyr}{ \!\!\!\!\!\!\!
    \left . \begin{array}{c} {} \\ {} \end{array}
    \right \} }
\nc{\longmid}{\left | \begin{array}{c} {} \\ {} \end{array}
    \right . \!\!\!\!\!\!\!}
\nc{\onetree}{\bullet}
\nc{\ora}[1]{\stackrel{#1}{\rar}}
\nc{\ola}[1]{\stackrel{#1}{\la}}
\nc{\ot}{\otimes}
\nc{\mot}{{{\boxtimes\,}}}
\nc{\otm}{\overline{\boxtimes}}
\nc{\sprod}{\bullet}
\nc{\scs}[1]{\scriptstyle{#1}}
\nc{\mrm}[1]{{\rm #1}}
\nc{\margin}[1]{\marginpar{\rm #1}}   
\nc{\dirlim}{\displaystyle{\lim_{\longrightarrow}}\,}
\nc{\invlim}{\displaystyle{\lim_{\longleftarrow}}\,}
\nc{\mvp}{\vspace{0.3cm}}
\nc{\tk}{^{(k)}}
\nc{\tp}{^\prime}
\nc{\ttp}{^{\prime\prime}}
\nc{\svp}{\vspace{2cm}}
\nc{\vp}{\vspace{8cm}}
\nc{\proofbegin}{\noindent{\bf Proof: }}
\nc{\proofend}{$\blacksquare$ \vspace{0.3cm}}
\nc{\modg}[1]{\!<\!\!{#1}\!\!>}
\nc{\intg}[1]{F_C(#1)}
\nc{\lmodg}{\!<\!\!}
\nc{\rmodg}{\!\!>\!}
\nc{\cpi}{\widehat{\Pi}}
\nc{\sha}{{\mbox{\cyr X}}}  
\nc{\shap}{{\mbox{\cyrs X}}} 
\nc{\shpr}{\diamond}    
\nc{\shp}{\ast}
\nc{\shplus}{\shpr^+}
\nc{\shprc}{\shpr_c}    
\nc{\msh}{\ast}
\nc{\zprod}{m_0}
\nc{\oprod}{m_1}
\nc{\vep}{\varepsilon}
\nc{\labs}{\mid\!}
\nc{\rabs}{\!\mid}

\nc{\mmbox}[1]{\mbox{\ #1\ }}
\nc{\fp}{\mrm{FP}} \nc{\rchar}{\mrm{char}} \nc{\End}{\mrm{End}} \nc{\Fil}{\mrm{Fil}}
\nc{\Mor}{Mor\xspace}
\nc{\gmzvs}{gMZV\xspace}
\nc{\gmzv}{gMZV\xspace}
\nc{\mzv}{MZV\xspace}
\nc{\mzvs}{MZVs\xspace}
\nc{\Hom}{\mrm{Hom}} \nc{\id}{\mrm{id}} \nc{\im}{\mrm{im}}
\nc{\incl}{\mrm{incl}} \nc{\map}{\mrm{Map}} \nc{\mchar}{\rm char}
\nc{\nz}{\rm NZ} \nc{\supp}{\mathrm Supp}

\nc{\Alg}{\mathbf{Alg}}
\nc{\Bax}{\mathbf{Bax}}
\nc{\bff}{\mathbf f}
\nc{\bfk}{{\bf k}}
\nc{\bfone}{{\bf 1}}
\nc{\bfx}{\mathbf x}
\nc{\bfy}{\mathbf y}
\nc{\base}[1]{\bfone^{\otimes ({#1}+1)}} 
\nc{\Cat}{\mathbf{Cat}}

\nc{\detail}{\marginpar{\bf More detail}
    \noindent{\bf Need more detail!}
    \svp}
\nc{\Int}{\mathbf{Int}}
\nc{\Mon}{\mathbf{Mon}}
\nc{\rbtm}{{shuffle }}
\nc{\rbto}{{Rota-Baxter }}
\nc{\remarks}{\noindent{\bf Remarks: }}
\nc{\Rings}{\mathbf{Rings}}
\nc{\Sets}{\mathbf{Sets}}

\nc{\BA}{{\Bbb A}} \nc{\CC}{{\Bbb C}} \nc{\DD}{{\Bbb D}}
\nc{\EE}{{\Bbb E}} \nc{\FF}{{\Bbb F}} \nc{\GG}{{\Bbb G}}
\nc{\HH}{{\Bbb H}} \nc{\LL}{{\Bbb L}} \nc{\NN}{{\Bbb N}}
\nc{\KK}{{\Bbb K}} \nc{\QQ}{{\Bbb Q}} \nc{\RR}{{\Bbb R}}
\nc{\TT}{{\Bbb T}} \nc{\VV}{{\Bbb V}} \nc{\ZZ}{{\Bbb Z}}

\nc{\cala}{{\mathcal A}} \nc{\calc}{{\mathcal C}}
\nc{\cald}{{\mathcal D}} \nc{\cale}{{\mathcal E}}
\nc{\calf}{{\mathcal F}} \nc{\calg}{{\mathcal G}}
\nc{\calh}{{\mathcal H}} \nc{\cali}{{\mathcal I}}
\nc{\call}{{\mathcal L}} \nc{\calm}{{\mathcal M}}
\nc{\caln}{{\mathcal N}} \nc{\calo}{{\mathcal O}}
\nc{\calp}{{\mathcal P}} \nc{\calr}{{\mathcal R}}
\nc{\cals}{{\mathcal S}}
\nc{\calt}{{\mathcal T}} \nc{\calw}{{\mathcal W}}
\nc{\calk}{{\mathcal K}} \nc{\calx}{{\mathcal X}}
\nc{\CA}{\mathcal{A}}

\nc{\fraka}{{\mathfrak a}}
\nc{\frakA}{{\mathfrak A}}
\nc{\frakb}{{\mathfrak b}}
\nc{\frakB}{{\mathfrak B}}
\nc{\frakC}{{\mathfrak C}}
\nc{\frakD}{{\mathfrak D}}
\nc{\frakg}{{\mathfrak g}}
\nc{\frakH}{{\mathfrak H}}
\nc{\frakL}{{\mathfrak L}}
\nc{\frakM}{{\mathfrak M}}
\nc{\bfrakM}{\overline{\frakM}}
\nc{\frakm}{{\mathfrak m}}
\nc{\frakP}{{\mathfrak P}}
\nc{\frakN}{{\mathfrak N}}
\nc{\frakp}{{\mathfrak p}}
\nc{\frakR}{{\mathfrak R}}
\nc{\frakS}{{\mathfrak S}}

\font\cyr=wncyr10
\font\cyrs=wncyr7

\begin{document}

\title[Distributive lattice, tilting and support $\tau$-tilting modules]
{Distributive lattices of tilting modules and support $\tau$-tilting modules over path algebras}

%
%
\author[Yichao Yang]{Yichao Yang}
\address{D\'{e}partement de math\'{e}matiques, Universit\'{e} de Sherbrooke, Sherbrooke, Qu\'{e}bec, Canada, J1K 2R1}
\email{yichao.yang@usherbrooke.ca}


\renewcommand{\thefootnote}{\alph{footnote}}
\setcounter{footnote}{-1} \footnote{}
\renewcommand{\thefootnote}{\alph{footnote}}
\setcounter{footnote}{-1} \footnote{\emph{2010 Mathematics Subject
Classification}: 16G20, 16G70, 05E10.}

\renewcommand{\thefootnote}{\alph{footnote}}
\setcounter{footnote}{-1} \footnote{\emph{Keywords}: Tilting module, $\tau$-tilting module, Distributive lattice, Auslander-Reiten quiver.}


\begin{abstract}
In this paper we study the poset of basic tilting $kQ$-modules when $Q$ is a Dynkin quiver, and the poset of basic support $\tau$-tilting $kQ$-modules when $Q$ is a connected acyclic quiver respectively. It is shown that the first poset is a distributive lattice if and only if $Q$ is of types $\mathbb{A}_{1}, \mathbb{A}_{2}$ or $\mathbb{A}_{3}$ with a nonlinear orientation and the second poset is a distributive lattice if and only if $Q$ is of type $\mathbb{A}_{1}$.
\end{abstract}

\maketitle

\setcounter{section}{0}

\section{Introduction}

Let $Q$ be a finite connected acyclic quiver and $kQ$ be the path algebra of $Q$ over an algebraically closed field $k$. Denote by {\rm mod}-$kQ$ the category of finite dimensional right $kQ$-modules, by {\rm ind}-$kQ$ the category of indecomposable modules in {\rm mod}-$kQ$ and by $\Gamma$({\rm mod} $kQ$) the Auslander-Reiten quiver of $kQ$. For $M\in$ {\rm mod}-$kQ$, we denote by {\rm add} $M$ (respectively, {\rm Fac} $M$, {\rm Sub} $M$) the category of all direct summands (respectively, factor modules, submodules) of finite direct sums of copies of $M$ and by $|M|$ the number of pairwise non-isomorphic indecomposable direct summands of $M$. Let $P_{i}$ be an indecomposable projective module in {\rm mod}-$kQ$ associated with vertex $i\in Q_{0}$ and $\tau$ be the Auslander-Reiten translation.

Tilting theory for $kQ$, or more generally for a finite dimensional basic $k$-algebra, was first appeared in \cite{[BB]} and have been central in the representation theory of finite dimensional algebras since the early seventies. For the classical tilting modules and their mutation theory, there is a naturally associated quiver named tilting quiver which is defined in \cite{[RS]}. Happel and Unger defined a partial order on the set of basic tilting modules and showed that the tilting quiver coincides with the Hasse quiver of this poset \cite{[HU]}. A related partial order has been studied in the $\tau$-tilting theory introduced in \cite{[AIR]} and the analog result also holds, that is, the support $\tau$-tilting quiver also coincide with the Hasse quiver of this related partial order.

Recently, the lattice structure of the poset of tilting modules and support $\tau$-tilting modules have been studied in \cite{[IRTT],[K],[R2]}. More precisely, Kase showed that for representation-infinite algebras $kQ$, the poset of its pre-projective tilting modules possess a distributive lattice structure if and only if the degree of all vertices in $Q$ are greater than $1$ \cite{[K]}. Later Iyama, Reiten, Thomas and Todorov proved that for path algebras $kQ$, the poset of its support $\tau$-tilting modules possess a lattice structure if and only if $Q$ is a Dynkin quiver or has at most $2$ vertices.

The aim of this paper is to study the following problem.

\begin{problem}\label{problem}
Let $Q$ be a finite connected acyclic quiver.

(1)~ When does the poset of basic tilting $kQ$-modules possess a distributive lattice structure?

(2)~ When does the poset of basic support $\tau$-tilting $kQ$-modules possess a distributive lattice structure?
\end{problem}

Our main result is the following theorem.

\begin{theorem}\label{main result}
Let $Q$ be a Dynkin quiver. Then the following statements are equivalent.

(1)~ All tilting modules are slice modules.

(2)~ The full subquiver generated by any tilting module form a section of $\Gamma$({\rm mod} $kQ$).

(3)~ The tilting quiver $\vec{\mathcal{T}}(Q)$ is a distributive lattice.

(4)~ Any boundary orbit (see Definition \ref{boundary}) of $\Gamma$({\rm mod} $kQ$) contains at most $2$ modules.
\end{theorem}

For the representation-infinite case, see \cite{[HV],[K],[KT]}.

As a consequence, the answer to Problem \ref{problem}(1) is given in the following theorem.

\begin{theorem}\label{answer to p1}
Let $Q$ be a finite connected acyclic quiver.

(1)~ {\rm [\cite{[K]}, Theorem 3.1]} If $Q$ is a non-Dynkin quiver, then the poset of basic pre-projective tilting $kQ$-modules is a distributive lattice if and only if the degree of all vertices in $Q$ are greater than $1$.

(2)~ If $Q$ is a Dynkin quiver, then the poset of basic tilting $kQ$-modules is a distributive lattice if and only if $Q$ is of types $\mathbb{A}_{1}, \mathbb{A}_{2}$ or $\mathbb{A}_{3}$ with a nonlinear orientation.
\end{theorem}

On the other hand, we also show the following result which answers Problem \ref{problem}(2).

\begin{theorem}\label{answer to p2}
Let $Q$ be a finite connected acyclic quiver. Then the poset of basic support $\tau$-tilting $kQ$-modules is a distributive lattice if and only if $Q$ is of type $\mathbb{A}_{1}$.
\end{theorem}

The paper is organized as follows. In section 2 we recall some preliminary definitions and results of tilting theory, $\tau$-tilting theory and lattice theory, especially about the tilting quiver, support $\tau$-tilting quiver and distributive lattice. In subsection 3.1 we first introduce the notions of boundary module and boundary orbit and then prove Theorem \ref{main result}. In subsection 3.2 we give a proof of Theorem \ref{answer to p1} by using Theorem \ref{main result}. In subsection 3.3 we prove Theorem \ref{answer to p2}.

\section{Preliminaries}

\subsection{Tilting theory and $\tau$-tilting theory}

We start with the following definitions of tilting modules and tilting quiver which was considered in \cite{[K]}, and was first introduced in \cite{[HU],[RS]}.

\begin{definition}
A module $T\in$ {\rm mod}-$kQ$ is a {\rm tilting module} if

(1)~ {\rm Ext}$_{kQ}^{1}(T,T)=0$.

(2)~ $|T|=|Q_{0}|$.
\end{definition}

We denote by $\mathcal{T}(Q)$ a complete set of representatives of the isomorphism classes of the basic tilting modules in {\rm mod}-$kQ$.

\begin{definition}
The {\rm tilting quiver} $\vec{\mathcal{T}}(Q)$ is defined as follows:

(1)~ $\vec{\mathcal{T}}(Q)_{0}:=\mathcal{T}(Q)$.

(2)~ $T\rightarrow T'$ in $\vec{\mathcal{T}}(Q)$ if $T\cong M\oplus X$, $T'\cong M\oplus Y$ for some $X,Y\in$ {\rm ind}-$kQ$, $M\in$ {\rm mod}-$kQ$ and there is a non-split exact sequence
$$\xymatrix{0 \ar[r] & X \ar[r] & M' \ar[r] & Y \ar[r] & 0}$$
with $M'\in$ {\rm add} $M$.
\end{definition}

Now we recall some basic definitions of $\tau$-tilting theory, which was first introduced in \cite{[AIR]}, in order to ``complete" the classical tilting theory from the viewpoint of mutation.

\begin{definition}
(1)~ We call $M\in$ {\rm mod}-$kQ$ {\rm $\tau$-rigid} if {\rm Hom}$_{kQ}(M,\tau M)=0$.

(2)~ We call $M\in$ {\rm mod}-$kQ$ {\rm $\tau$-tilting} if $M$ is $\tau$-rigid and $|M|=|Q_{0}|$.

(3)~ We call $M\in$ {\rm mod}-$kQ$ {\rm support $\tau$-tilting} if there exists an idempotent $e$ of $kQ$ such that $M$ is a $\tau$-tilting $(kQ/ \langle e\rangle)$-module.
\end{definition}

We denote by $\mathcal{ST}(Q)$ a complete set of representatives of the isomorphism classes of the basic support $\tau$-tilting modules in {\rm mod}-$kQ$.

Recall that the Hasse-quiver $\vec{P}$ of a poset $(P,\leq)$ is defined as follows:

(1)~ $\vec{P}_{0}:=P$.

(2)~ $x\rightarrow y$ in $\vec{P}$ if $x>y$ and there is no $z\in P$ such that $x>z>y$.

The support $\tau$-tilting quiver $\vec{\mathcal{ST}}(Q)$ is defined as follows.

\begin{prop-def}[\cite{[AIR]}, Theorem 2.7, Corollary 2.34]\label{Hasse1}
(1)~ Let $T,T'\in \mathcal{ST}(Q)$, then the following relation $\leq$ defines a partial order on $\mathcal{ST}(Q)$,
$$T\geq T' \stackrel{\rm def}{\Leftrightarrow} {\rm Fac} T\supseteq {\rm Fac} T'$$.

(2)~ The {\rm support $\tau$-tilting quiver} $\vec{\mathcal{ST}}(Q)$ is the Hasse quiver of the partial order set $(\mathcal{ST}(Q),\leq)$.
\end{prop-def}

We remark that there is the following similar result in the classical tilting theory.

\begin{theorem}[\cite{[HU]}, Theorem 2.1]\label{Hasse2}
(1)~ Let $T,T'\in \mathcal{T}(Q)$, then the following relation $\leq$ defines a partial order on $\mathcal{T}(Q)$,
$$T\geq T' \stackrel{\rm def}{\Leftrightarrow} {\rm Fac} T\supseteq {\rm Fac} T'$$.

(2)~ The tilting quiver $\vec{\mathcal{T}}(Q)$ is the Hasse quiver of the partial order set $(\mathcal{T}(Q),\leq)$.
\end{theorem}

We end this subsection with the following two examples.

\begin{example}\label{ex1}
Let $Q_{1}, Q_{2}$ be the following two different quivers, see Figure 1. Although they share the same underlying graph, however, the corresponding tilting quivers are different.

\begin{figure}[h] \centering

  \includegraphics*[83,565][481,769]{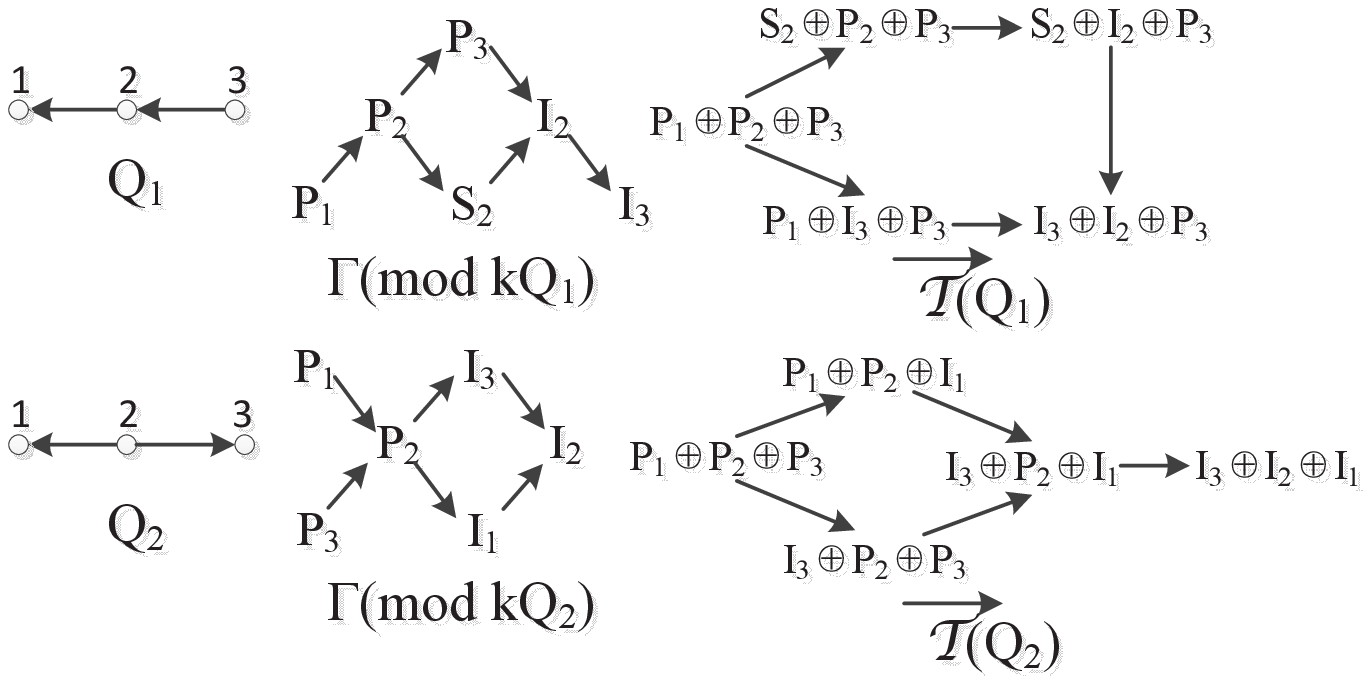}

  {\rm Figure 1}
\end{figure}
\end{example}

\begin{example}\label{ex2}
Let $Q$ be of type $\mathbb{A}_{2}$, then its support $\tau$-tilting quiver $\vec{\mathcal{ST}}(Q)$ is shown in Figure 2.

\begin{figure}[h] \centering

  \includegraphics*[85,666][333,765]{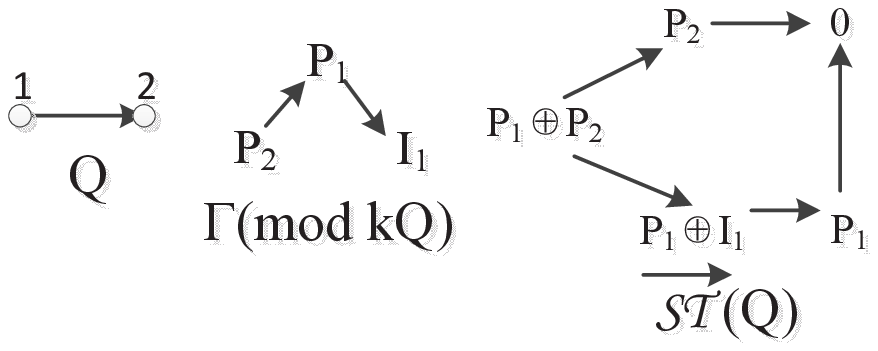}

  {\rm Figure 2}
\end{figure}
\end{example}

\subsection{Lattices and distributive lattices}

In this subsection we will recall definitions of lattices and distributive lattices.

\begin{definition}
A poset $(L,\leq)$ is a {\rm lattice} if for any $x,y\in L$ there is a minimum element of $\{z\in L | z\geq x,y\}$ and there is a maximum element of $\{z\in L | z\leq x,y\}$.

In this case, we denote by $x\vee y$ the minimum element of $\{z\in L | z\geq x,y\}$ and call it {\rm join} of $x$ and $y$. We also denote by $x\wedge y$ the maximum element of $\{z\in L | z\leq x,y\}$ and call it {\rm meet} of $x$ and $y$.
\end{definition}

\begin{definition}
A lattice $L$ is a {\rm distributive lattice} if $(x\vee y)\wedge z=(x\wedge z)\vee (y\wedge z)$ holds for any $x,y,z\in L$.
\end{definition}

Immediately we have the following basic observation, which will be used frequently in this paper.

\begin{lemma}\label{easy lemma}
For any $n\geq 2$, the following Hasse quiver in Figure 3 is not a distributive lattice.

\begin{figure}[h] \centering

  \includegraphics*[104,707][273,776]{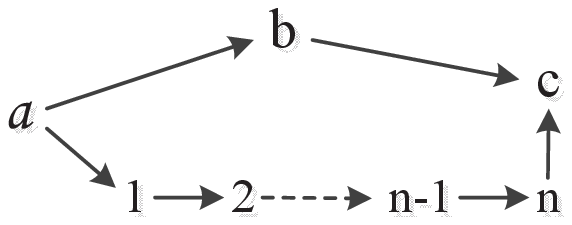}

  {\rm Figure 3}
\end{figure}
\end{lemma}

\begin{proof}
Since $n\geq 2$, it is easy to see that
$$(b\vee 2)\wedge 1=a\wedge 1=1\neq 2=c \vee 2=(b\wedge 1)\vee (2 \wedge 1),$$
therefore it is not a distributive lattice.
\end{proof}

In the above examples \ref{ex1} and \ref{ex2}, it is easy to see that the lattice $(\mathcal{T}(Q_{2}),\leq)$ is a distributive lattice. On the other hand, it follows by Lemma \ref{easy lemma} that both $(\mathcal{T}(Q_{1}),\leq)$ and $(\vec{\mathcal{ST}}(Q),\leq)$ are not distributive lattice.

\section{Main results}

\subsection{Boundary module and boundary orbit.}

From now on, we will not distinguish between an indecomposable $kQ$-module $M$ and its corresponding vertex $[M]$ in the Auslander-Reiten quiver $\Gamma$({\rm mod} $kQ$). We will also not distinguish between a poset $(P,\leq)$ and its Hasse quiver $\vec{P}$.

By Theorem \ref{Hasse2} and Proposition-Definition \ref{Hasse1}, it is easy to see that our problem reduces to the study of lattice structure of the tilting quiver $\vec{\mathcal{T}}(Q)$ and the support $\tau$-tiling quiver $\vec{\mathcal{ST}}(Q)$.

Before proceeding further, let $(\Gamma,\tau)$ be a connected translation quiver, recall from \cite{[ABS]} that a connected full subquiver $\Sigma$ of $\Gamma$ is called a \emph{presection} (is also called a \emph{cut} in \cite{[L2]}) in $\Gamma$ if it satisfies the following two conditions:

(1)~ If $x\in \Sigma_{0}$ and $x\rightarrow y$ is an arrow, then either $y\in \Sigma_{0}$ or $\tau y\in \Sigma_{0}$.

(2)~ If $y\in \Sigma_{0}$ and $x\rightarrow y$ is an arrow, then either $x\in \Sigma_{0}$ or $\tau^{-1}x\in \Sigma_{0}$.

Moreover, in \cite{[L1]} a connected full subquiver $\Sigma$ of $\Gamma$ is a called \emph{section} of $\Gamma$ if the following conditions are satisfied:

(1)~ $\Sigma$ contains no oriented cycle.

(2)~ $\Sigma$ meets each $\tau$-orbit in $\Gamma$ exactly once.

(3)~ $\Sigma$ is convex in $\Gamma$, that is, every path in $\Gamma$ with end-points belonging to $\Sigma$ lies entirely in $\Sigma$.

From \cite{[R1]} recall also that a module $S$ is said to be a \emph{slice module} if $S$ is sincere and {\rm add} $S$ satisfies the following conditions:

(1)~ If there is a path $x_{0}\rightarrow x_{1}\rightarrow \cdots \rightarrow x_{t}$ with $x_{0},x_{t}\in$ {\rm add} $S$ in the Auslander-Reiten quiver, then $x_{i}\in$ {\rm add} $S$ ($i=0,1,\cdots,t$).

(2)~ If $M$ is indecomposable and not projective, then at most one of $M, \tau M$ belongs to {\rm add} $S$.

(3)~ If there is an arrow $M\rightarrow X$ with $X\in$ {\rm add} $S$ in the Auslander-Reiten quiver, then either $M\in$ {\rm add} $S$ or $M$ is not injective and $\tau^{-1} M\in$ {\rm add} $S$.

Now we introduce the notions of boundary module and boundary orbit.

\begin{definition}\label{boundary}
(1)~ We call a module $M\in \Gamma$({\rm mod} $kQ$) {\rm boundary module} if $M$ has at most one direct predecessor and at most one direct successor in Auslander-Reiten quiver $\Gamma$({\rm mod} $kQ$).

(2)~ We call a $\tau$-orbit $\Sigma$ of $\Gamma$({\rm mod} $kQ$) {\rm boundary orbit} if $\Sigma$ contains a boundary module.
\end{definition}

The following observation is useful.

\begin{lemma}\label{useful lemma}
Let $Q$ be a Dynkin quiver. If one of its boundary orbits contains at least $3$ modules, then the tilting quiver $\vec{\mathcal{T}}(Q)$ is not a distributive lattice.
\end{lemma}

\begin{proof}
Since $Q$ is a Dynkin quiver, $\Gamma$({\rm mod} $kQ$) must be a full convex subquiver of $\mathbb{Z}Q$. Without loss of generality, by our assumption $\Gamma$({\rm mod} $kQ$) will contain the following shaded area $\mathcal{T}$, see Figure 4.

\begin{figure}[h] \centering

  \includegraphics*[11,508][445,778]{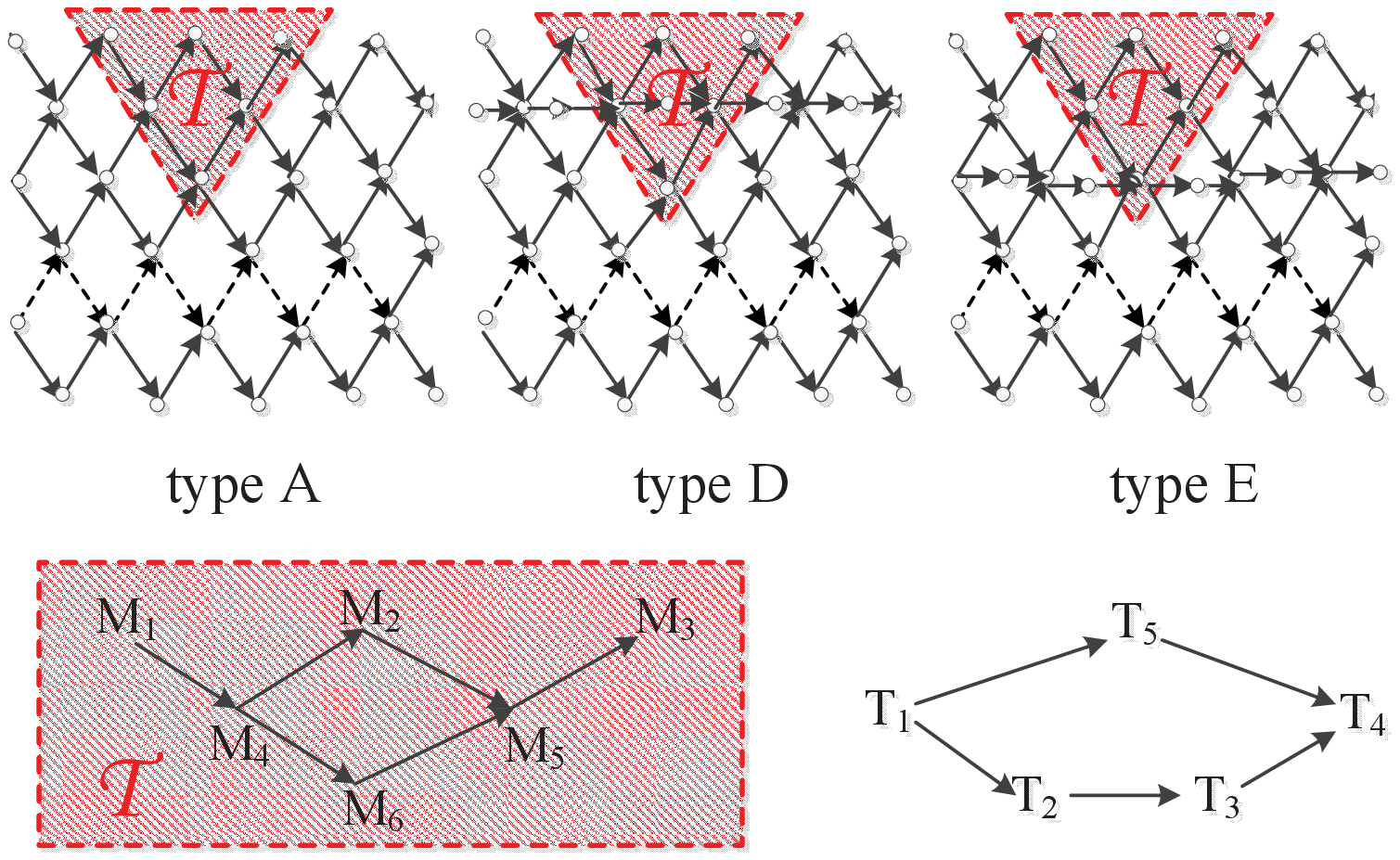}

  {\rm Figure 4}
\end{figure}

Now we enlarge $\mathcal{T}$ for each type, for the type $A$, see the left-lower of Figure 4. For simplicity, we may continue with the type $A$, for the remaining two types, the argument is similar.

Let $|Q_{0}|=n$, it is easy to see that we can construct a section $\Sigma$ of the lower $(n-2)$-rows starting with $M_{6}$ and denote the module corresponding to this section by $M_{\Sigma}$. Then we consider the following five modules
$$T_{1}=M_{\Sigma}\oplus M_{4}\oplus M_{1}, T_{2}=M_{\Sigma}\oplus M_{4}\oplus M_{2}, T_{3}=M_{\Sigma}\oplus M_{5}\oplus M_{2},$$
$$T_{4}=M_{\Sigma}\oplus M_{5}\oplus M_{3}, T_{5}=M_{\Sigma}\oplus M_{1}\oplus M_{3}.$$

Since $\Gamma$({\rm mod} $kQ$) is a standard component, it is not hard to see that all of these five modules are tilting modules and they forms the right-lower of Figure 4, which is a full subquiver of the tilting quiver $\vec{\mathcal{T}}(Q)$, however, is not a distributive lattice by Lemma \ref{easy lemma}. Hence the tilting quiver $\vec{\mathcal{T}}(Q)$ is also not a distributive lattice, which completes the proof.
\end{proof}

Now we are ready to prove Theorem \ref{main result}.

(1) $\Leftrightarrow$ (2): This is shown in \cite{[LY]} or \cite{[Z]}.

(2) $\Rightarrow$ (3): Let $|Q_{0}|=n$, according to (2) it follows that any tilting module can be written as
$$T\cong \bigoplus_{i=1}^{n}\tau^{-r_{i}} P_{i}$$
for $r_{i}\in \mathbb{Z}_{\geq 0}, 1\leq i\leq n$ and if $T,T'$ be two tilting modules, $T\rightarrow T'$ in $\vec{\mathcal{T}}(Q)$ if and only if there is an indecomposable direct summand $X$ such that $T\cong M\oplus X$ and $T'\cong M\oplus \tau^{-1}X$. Thus, for any two tilting modules $T\cong \bigoplus_{i=1}^{n}\tau^{-r_{i}} P_{i}, T'\cong \bigoplus_{i=1}^{n}\tau^{-r_{i}'} P_{i}$, $T\geq T'$ if and only if $r_{i}\leq r_{i}', 1\leq i\leq n$.

From now on let $\Sigma_{T}$ be the full subquiver of $\Gamma$({\rm mod} $kQ$) generated by $T$. Since $\Sigma_{T},\Sigma_{T'}$ form a section of $\Gamma$({\rm mod} $kQ$), it is not hard to check that both $\Sigma_{\bigoplus_{i=1}^{n}\tau^{-{\rm min}\{r_{i},r_{i}'\}} P_{i}}$ and $\Sigma_{\bigoplus_{i=1}^{n}\tau^{-{\rm max}\{r_{i},r_{i}'\}} P_{i}}$ again form a section of $\Gamma$({\rm mod} $kQ$), which implies that both $\bigoplus_{i=1}^{n}\tau^{-{\rm min}\{r_{i},r_{i}'\}} P_{i}$ and $\bigoplus_{i=1}^{n}\tau^{-{\rm max}\{r_{i},r_{i}'\}} P_{i}$ are tilting modules. Therefore the join and meet of $T$ and $T'$ are
$$T\vee T'\cong \bigoplus_{i=1}^{n}\tau^{-{\rm min}\{r_{i},r_{i}'\}} P_{i},~~~~~~T\wedge T'\cong \bigoplus_{i=1}^{n}\tau^{-{\rm max}\{r_{i},r_{i}'\}} P_{i}$$
respectively, which makes the tilting quiver $\vec{\mathcal{T}}(Q)$ to be a distributive lattice. Indeed, it follows by the fact that $a\vee b=({\rm min}(r_{i},r_{i}'))_{1\leq i\leq n}$ and $a\wedge b=({\rm max}(r_{i},r_{i}'))_{1\leq i\leq n}$ makes $(\mathbb{Z}^{n},\leq^{{\rm op}})$ to be a distributive lattice, where $a=(r_{i})_{1\leq i\leq n}, b=(r_{i}')_{1\leq i\leq n}$.

(3) $\Rightarrow$ (4): It follows from Lemma \ref{useful lemma} at once.

(4) $\Rightarrow$ (2): Since any boundary orbit of $\Gamma$({\rm mod} $kQ$) contains at most $2$ modules and $\Gamma$({\rm mod} $kQ$) is a full convex subquiver of $\mathbb{Z}Q$, it follows that $\Gamma$({\rm mod} $kQ$) is bounded by the following shaded area $\mathcal{R}$, see Figure 5.

\begin{figure}[h] \centering

  \includegraphics*[28,622][500,788]{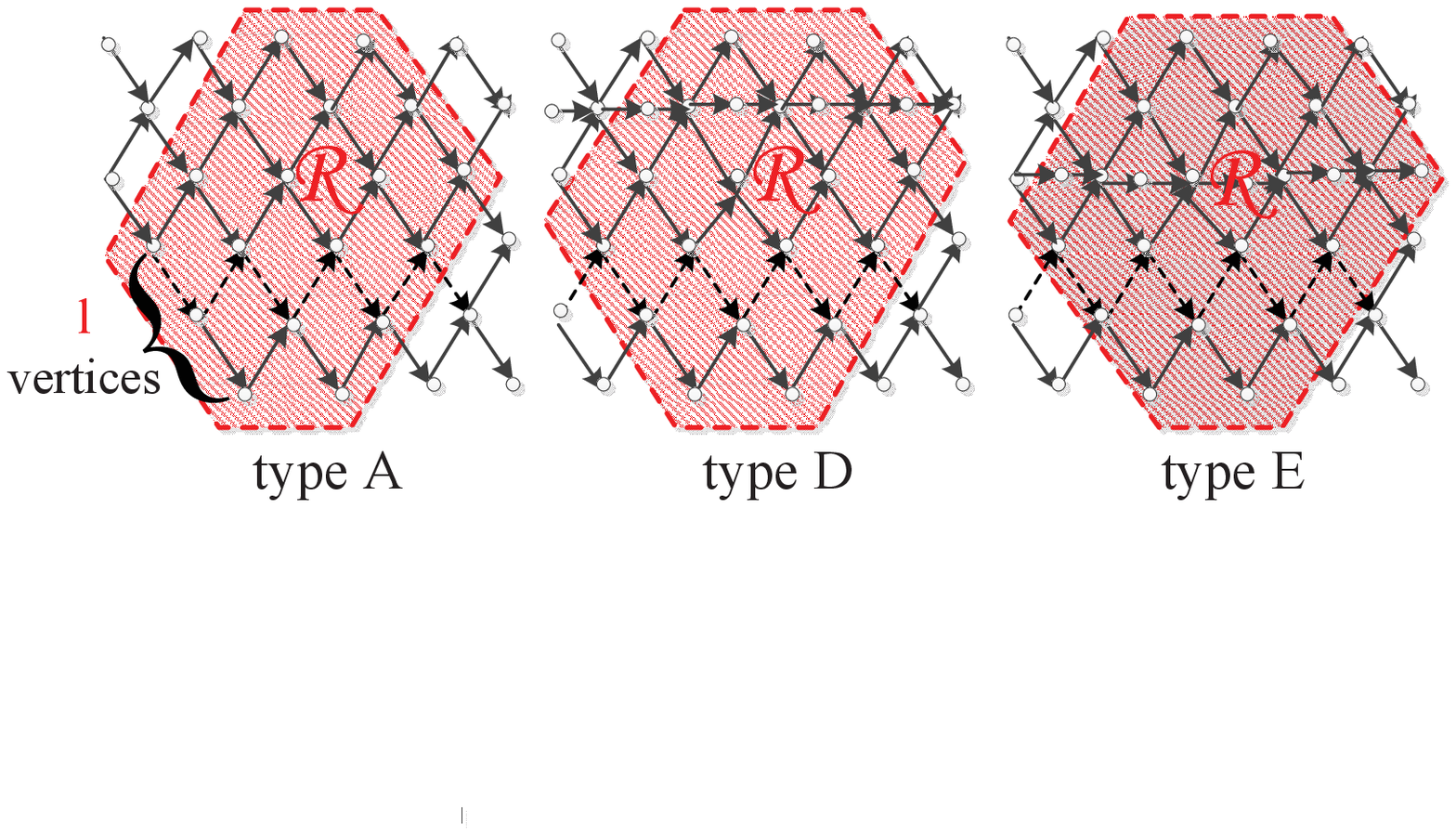}

  {\rm Figure 5}
\end{figure}

Since $\Gamma$({\rm mod} $kQ$) is a standard component, we have that for any $M,N\in \mathcal{R}$, if there exists a path from $M$ to $\tau N$, then {\rm Hom}$_{kQ}(M,\tau N)\neq 0$.

Let $T$ be any tilting module, because {\rm Ext}$_{kQ}^{1}(T,T)=${\rm Hom}$_{kQ}(T,\tau T)=0$, so there is no path from $T_{i}$ to $\tau T_{j}$, which implies that $\Sigma_{T}$ meets each $\tau$-orbit at most once. Moreover, since $|(\Sigma_{T})_{0}|=|T|=|Q_{0}|$, it follows that $\Sigma_{T}$ meets each $\tau$-orbit exactly once.

According to [\cite{[ABS]}, Proposition 1.7], it suffices to prove that $\Sigma_{T}$ is a presection of $\Gamma$({\rm mod} $kQ$). Indeed, if $x\in (\Sigma_{T})_{0}$, $x\rightarrow y$ is an arrow and $y, \tau y\notin (\Sigma_{T})_{0}$, then there exists $i\neq 0,1$ such that $\tau^{i} y\in (\Sigma_{T})_{0}$. If $i\geq 2$, then there exists a path from $\tau^{i}y$ to $\tau x$, $\tau^{i}y, x\in (\Sigma_{T})_{0}\subseteq \Gamma$({\rm mod} $kQ$)$\subseteq \mathcal{R}$, thus we have {\rm Hom}$_{kQ}(\tau^{i}y,\tau x)\cong$ {\rm Ext}$_{kQ}^{1}(x,\tau^{i}y)\neq 0$, which contradicts that $\tau^{i}y, x\in (\Sigma_{T})_{0}$ and $T$ is a tilting module. For the $i\leq -1$ case, the proof is similar.

Using the same argument as above, we can easily carry out that if $y\in (\Sigma_{T})_{0}$, $x\rightarrow y$ is an arrow, then either $x\in (\Sigma_{T})_{0}$ or $\tau^{-1}x\in (\Sigma_{T})_{0}$. Finally the connectivity of $\Sigma_{T}$ follows from the connectivity of $\Gamma$({\rm mod} $kQ$), which completes the proof.

\subsection{Proof of Theorem \ref{answer to p1}}

In this subsection we start to prove Theorem \ref{answer to p1}.

For the non-Dynkin case, see [\cite{[K]}, Theorem 3.1]. If $Q$ is a Dynkin quiver, we divide into the following three cases.

{\bf Case 1: $Q$ is of type $A$.}

$|Q_{0}|=1,2$, then the tilting quivers are $\cdot,\cdot\rightarrow\cdot$, respectively, it is clear.

$|Q_{0}|=3$, see Example \ref{ex1} and it is easy to see that the tilting quiver of $\cdot\rightarrow\cdot\leftarrow\cdot$ is also a distributive lattice.

$|Q_{0}|=4$, then we can list all the non-isomorphic quivers and their corresponding Auslander-Reiten quivers as follows, see Figure 6.

\begin{figure}[h] \centering

  \includegraphics*[47,600][340,760]{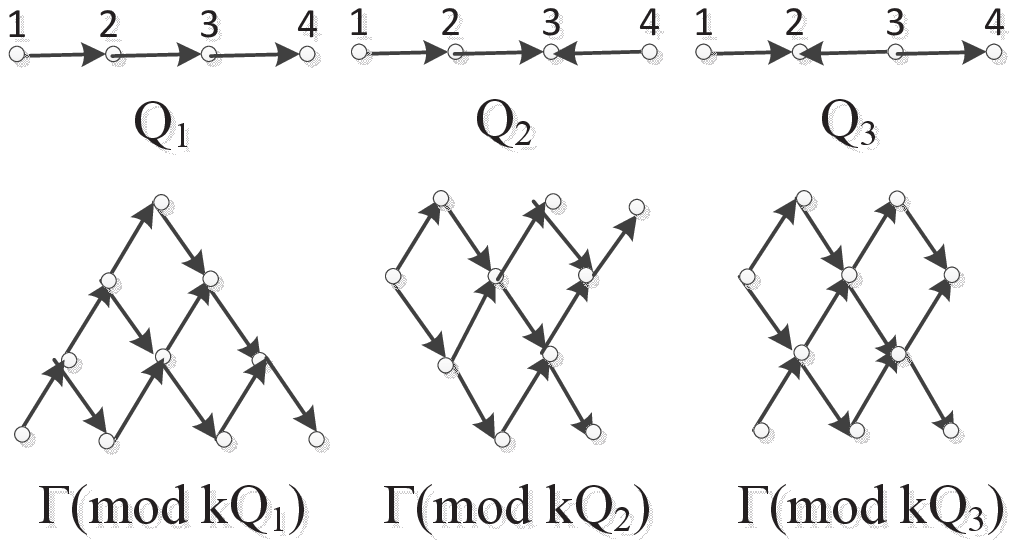}

  {\rm Figure 6}
\end{figure}

Since for each of these three Auslander-Reiten quivers, we can always find a boundary orbit containing $3$ modules, then by Theorem \ref{main result} the corresponding tilting quiver $\vec{\mathcal{T}}(Q_{i})$ is not a distributive lattice, $1\leq i\leq 3$.

$|Q_{0}|\geq 5$, if the tilting quiver $\vec{\mathcal{T}}(Q)$ is a distributive lattice, then by Theorem \ref{main result} any boundary orbit of $\Gamma$({\rm mod} $kQ$) contains at most $2$ modules, i.e., $\Gamma$({\rm mod} $kQ$) is bounded by the shaded area $\mathcal{R}$ of Figure 5.

Let $|Q_{0}|=n\geq 5$ and $l$ be defined in Figure 5. It is well known that the number of indecomposable $kQ$-modules is $\frac{n(n+1)}{2}$. On the other hand, there are at most $l(n+1-l)+n$ modules in $\mathcal{R}$, $1\leq l\leq n$. However, when $n\geq 5$ we have
$$l(n+1-l)+n=-(l-\frac{n+1}{2})^{2}+\frac{n^{2}+6n+1}{4}\leq \frac{n^{2}+6n+1}{4}<\frac{n(n+1)}{2}$$
which contradicts that $\Gamma$({\rm mod} $kQ$) is bounded by $\mathcal{R}$.

{\bf Case 2: $Q$ is of type $D$.}

Similarly, if the tilting quiver $\vec{\mathcal{T}}(Q)$ is a distributive lattice, then $\Gamma$({\rm mod} $kQ$) is bounded by $\mathcal{R}$. Let $|Q_{0}|=n\geq 4$ and $l$ is defined in the same way, then on one hand the number of indecomposable $kQ$-modules is $n(n-1)$; On the other hand, there are at most $l(n-l)+n+3$ modules in $\mathcal{R}$, $1\leq l\leq n-1$. However, when $n\geq 4$ we have
$$l(n-l)+n+3=-(l-\frac{n}{2})^{2}+\frac{n^{2}+4n+12}{4}\leq \frac{n^{2}+4n+12}{4}<n(n-1)$$
the same contradiction follows.

{\bf Case 3: $Q$ is of type $E$.}

We now proceed as in the proof of above two cases. On one hand, when $n=6,7,8$, the number of indecomposable $kQ$-modules is $36,63,120$, respectively. On the other hand, there are at most $l(n-l)+n+4$ modules in $\mathcal{R}$, $1\leq l\leq n-1$. However,
$$l(n-l)+n+4=-(l-\frac{n}{2})^{2}+\frac{n^{2}+4n+16}{4}\leq \frac{n^{2}+4n+16}{4}$$
which equals to $19,23.25,28$ respectively when $n=6,7,8$, now we have the same contradiction.

Finally, by combining the above three cases together, we complete the proof of Theorem \ref{answer to p1}(2).

\subsection{Proof of Theorem \ref{answer to p2}}

Indeed, by [\cite{[IRTT]}, Theorem 0.4] it suffices to consider the following two cases.

{\bf Case 1: $Q$ is of Dynkin type.}

If $|Q_{0}|=1$, then the support $\tau$-tilting quiver is $\cdot\rightarrow\cdot$, it is clear.

If $|Q_{0}|=n\geq 2$, then $Q$ contains $\mathbb{A}_{2}$ as its full subquiver. Without loss of generality we assume that $\{e_{1},\cdots,e_{n}\}$ is a complete set of primitive orthogonal idempotents for $kQ$ and there is an arrow $\alpha$ between the vertices $1$ and $2$. Let $e=e_{3}+e_{4}+\cdots+e_{n}$, then $kQ/\langle e\rangle\cong k\mathbb{A}_{2}$.

By Example \ref{ex2} the support $\tau$-tilting quiver $\vec{\mathcal{ST}}(\mathbb{A}_{2})$ is not a distributive lattice. On the other hand, according to [\cite{[AIR]}, Proposition 2.27(a)] it can easily be seen that $\vec{\mathcal{ST}}(\mathbb{A}_{2})$ is a full subquiver of $\vec{\mathcal{ST}}(Q)$, which implies that $\vec{\mathcal{ST}}(Q)$ is not a distributive lattice itself.

{\bf Case 2: $Q$ has at most $2$ vertices.}

According to [\cite{[IRTT]}, Proposition 2.2], it follows that the support $\tau$-tilting quiver $\vec{\mathcal{ST}}(Q)$ is isomorphic to the Figure 3 in Lemma \ref{easy lemma}, where $n$ tends to $+\infty$. Now by Lemma \ref{easy lemma} it is obvious that $\vec{\mathcal{ST}}(Q)$ is not a distributive lattice.

Finally, by combining the above two cases together, we complete the proof of Theorem \ref{answer to p2}.

\bigskip

{\bf Acknowledgements.}~~ This work was carried out when the author is a postdoctoral fellow at Universit\'{e} de Sherbrooke, financed by Fonds Qu\'{e}b\'{e}cois de la Recherche sur la Nature et les Technologies (Qu\'{e}bec, Canada) through the Merit Scholarship Program For Foreign Students. He would like to thank Professor Shiping Liu for his valuable discussions.

\end{document}